\documentclass[12pt]{amsart}
\usepackage{latexsym,amsfonts,amsthm,amsmath,amscd,amssymb}
\usepackage[dvips]{graphicx}

\def\Ker{\mbox{\rm Ker}}

\newcommand{\Z} {{\mathbb Z}}
\newcommand{\Hy} {\mathbb H}
\newcommand{\R} {\mathbb R}
\newcommand{\matP} {\mathbb P}
\newcommand{\calL}{\mathcal L}

\def\emline#1#2#3#4#5#6{
       \put(#1,#2){\special{em:moveto}}
       \put(#4,#5){\special{em:lineto}}}
\newtheorem{lemma}{Lemma}[section]
\newtheorem{thm}[lemma]{Theorem}

\newtheorem{prop}[lemma]{Proposition}
\newtheorem{cor}[lemma]{Corollary}

\begin{document}

\author[Matveev]{Sergei Matveev}
\address{Chelyabinsk State University, Chelyabinsk 454021, Russia}
\email{matveev@csu.ru}

\author[Petronio]{Carlo Petronio}
\address{Departimento di Matematica Applicata, Universit\`a di Pisa,
Via Bonanno Pisano 25B,
56126, Pisa, Italy} \email{petronio@dm.unipi.it}

\author[Vesnin]{Andrei Vesnin}
\address{Sobolev Institute of Mathematics, Novosibirsk 630090, Russia}
\email{vesnin@math.nsc.ru}

\title[Bounds for the complexity of hyperbolic 3-manifolds]
{Two-sided asymptotic bounds for the complexity of some closed
hyperbolic three-manifolds}

\thanks{This work is the result of a collaboration among the three
authors carried out in the frame of the INTAS project
``CalcoMet-GT'' 03-51-3663. The first named author was also
supported by the Russian Fund for Fundamental Research, grant
05-01-00293.}



\begin{abstract}
We establish two-sided bounds for the complexity of two infinite
series of closed orientable 3-dimensional hyperbolic manifolds, the
L\"obell manifolds and the Fibonacci manifolds.
\end{abstract}

\date{\today}

\maketitle

\section{Introduction}
If $M$ is a compact 3-dimensional manifold, its
\emph{complexity}~\cite{DAN, Acta} is a non-negative integer $c(M)$
which formally translates the intuitive notion of ``how
complicated'' $M$ is. In particular, if $M$ is closed and
irreducible and different from the 3-sphere $\mathbb S^3$, the
projective 3-space $\R\matP^3$, and the lens space $L(3,1)$, its
complexity $c(M)$ is equal to the minimum of the number of
tetrahedra over all ``singular'' triangulations of $M$. (A singular
triangulation of $M$ is a realization of $M$ as a union of
tetrahedra with pairwise glued 2-faces). The complexity function has
many natural properties, among which additivity under connected sum.

The task of computing the complexity $c(M)$ of a given manifold $M$
is extremely difficult. For closed $M$, the exact value is presently
known only if $M$ belongs to the computer-generated tables of
manifolds up to complexity 12, see~\cite{DAN05}. Therefore the
problem of finding ``reasonably good'' two-sided bounds for $c(M)$
is of primary importance. The first results of this kind were
obtained in~\cite{mape,pepe}, where an estimate on $c(M)$ was given
in terms of the properties of the homology groups of~$M$.

In the present paper we establish two-sided bounds on the complexity
of two infinite series of closed orientable 3-dimensional hyperbolic
manifolds, the L\"obell manifolds and the Fibonacci manifolds. The
upper bounds are obtained by constructing fundamental polytopes of
these manifolds in hyperbolic space $\Hy^3$, while the lower bounds
(which are only proved in an ``asymptotic'' fashion, see below) are
based on the calculation of their volumes. We mention here that the
L\"obell manifolds are constructed from polytopes that generalize
the right-angled dodecahedron, and the Fibonacci manifolds are
constructed from polytopes that generalize the regular icosahedron.

Before turning to the statements and proofs of our estimates, we
remark that, in the class of compact 3-manifolds with non-empty
boundary, exact values of complexity are currently known for two
infinite families. The first one consists of the manifolds which
finitely cover the complement of the figure-eight knot or its
``twin'' (a different manifold with the same
volume)~\cite{anis,matbook}. The second family consists of the
manifolds $M$ such that $\partial M$ consists of $k\geqslant0$ tori
and a surface of genus $g\geqslant2$, and $M$ admits an ideal
triangulation with $g+k$ tetrahedra~\cite{fmp1,fmp2}. 

The authors
thank Ekaterina Pervova for her help.

\section{L\"obell manifolds}

In this section we obtain upper and lower bounds on manifold
complexity for a certain infinite family of closed hyperbolic
3-manifolds which generalize the classical L\"obell manifold. Recall
that in order to give a positive answer to the question of the
existence of ``Clifford-Klein space forms'' (that is, closed
manifolds) of constant negative curvature, F.~L\"obell~\cite{Lo}
constructed in 1931 the first example of a closed orientable
hyperbolic 3-manifold. The manifold was obtained by gluing together
eight copies of the right-angled 14-faced polytope (denoted by
$R(6)$ and shown in Fig.~\ref{Fig1} below) with an upper and a lower
basis both being regular hexagons, and a lateral surface given by 12
pentagons, arranged similarly as in the dodecahedron.

An algebraic approach to constructing hyperbolic 3-manifolds from
eight copies of a right-angled polytope was suggested in~\cite{Ve1},
as we will now describe. Let us fix in hyperbolic 3-space $\Hy^3$ a
bounded polytope $R$, namely a compact and convex set homeomorphic
to the 3-disc, with boundary given by a finite union of geodesic
polygons, called the faces of $R$. The notions of vertex and edge of
$R$ are defined in the obvious fashion. We further assume that $R$
is right-angled, namely that all the dihedral angles along the edges
of $R$ are $\pi/2$, which easily implies that all the faces of $R$
have at least 5 edges and all vertices of $R$ are trivalent. Note
that, according to Andreev's theorem~\cite{Andre}, these
combinatorial conditions on an abstract polytope are actually also
sufficient for its realizability as a bounded right-angled polytope
in $\Hy^3$.

We will denote henceforth by $G$ the subgroup of ${\rm
Isom}(\Hy^3)$, the isometry group of hyperbolic 3-space, generated
by the reflections in the planes containing the faces of $R$. The
following is an easy consequence of Poincar\'e's polyhedron
theorem~\cite{EpPet,Vin}:

\begin{lemma}\label{Rfund:lem}
$R$ is a fundamental domain for $G$ and a presentation of $G$ is given by:
\begin{itemize}
\item A generator for each face of $R$;
\item The relation $\rho^2$ for each
generator $\rho$ and the relation $[\rho_1,\rho_2]$ for each pair of
generators $\rho_1,\rho_2$ associated to faces sharing an edge.
\end{itemize}
\end{lemma}

This result implies that $G$ is a discrete subgroup of isometries of
$\Hy^3$. In particular, a subgroup $K$ of $G$ will act freely on
$\Hy^3$ if and only if it is torsion-free. Moreover for each vertex
$v$ of $R$ the stabilizer $\mbox{Stab}_G (v)$ of $v$ in $G$ is
isomorphic to the Abelian group $\Z/_{\!2} \oplus \Z/_{\!2} \oplus
\Z/_{\!2} = (\Z/_{\!2})^3$ of order $8$, which we will view as a
vector space over the field $\Z/_{\!2}$. We next quote two lemmas
proved in~\cite{Ve1} and derive an easy consequence.

\begin{lemma}  \label{lemmaL2.1}
If $\varphi : G \rightarrow (\Z/_{\!2})^3$ is an epimorphism, the
following are equivalent:
\begin{enumerate}
\item $\Ker(\varphi)$ is torsion-free; \item For each vertex of
$R$, if $\rho_1,\rho_2,\rho_3$ are the reflections in the faces of
$R$ incident to the vertex,
$\varphi(\rho_1),\varphi(\rho_2),\varphi(\rho_3)$ are linearly
independent over $\Z/_{\!2}$.\end{enumerate}
\end{lemma}

We consider now in $(\Z/_{\!2})^3$ the vectors $\alpha = (1, 0, 0)$,
$\beta = (0, 1, 0)$, $\gamma = (0, 0, 1)$ and $\delta = \alpha +
\beta + \gamma = (1, 1, 1)$, and we note that any three of them are
linearly independent.

\begin{lemma}  \label{lemmaL2.2}
Let $\varphi : G \rightarrow (\Z/_{\!2})^3$ be an epimorphism that
maps each of the generating reflections of $G$ to one of the
elements $\alpha$, $\beta$, $\gamma$, $\delta$. Then $\Ker(\varphi)$
is a subgroup of ${\rm Isom}^+(\Hy^3)$, the group of
orientation-preserving isometries of hyperbolic 3-space.
\end{lemma}

\begin{prop} \label{easy:prop} If an epimorphism $\varphi: G \rightarrow
(\Z/_{\!2})^3$ satisfies condition (2) of Lemma~\ref{lemmaL2.1} and
the hypothesis of~\ref{lemmaL2.2} then the quotient $M = \Hy^3 /
\Ker(\varphi)$ is a closed orientable hyperbolic $3$-manifold.
\end{prop}

\begin{proof}
Lemma~\ref{lemmaL2.1} implies that $\Ker(\varphi)$ is torsion-free,
so the quotient is a hyperbolic 3-manifold $M$ without boundary, and
Lemma~\ref{lemmaL2.2} implies that $M$ is orientable. Since
$\Ker(\varphi)$ has index $8$ in $G$ and $R$ is a fundamental domain
for $G$, a fundamental domain for $\Ker(\varphi)$ is given by
$\bigcup_{i=1}^8 g_i(R)$ where $\{g_i:\ 1\leqslant i\leqslant 8\}$ is a set of
representatives of $G/\Ker(\varphi)$. Such a fundamental domain is
compact, so $M$ is compact, whence closed.
\end{proof}

If we now describe a homomorphism of $G$ by labelling each face of
$R$ by the image of the reflection in the plane containing that
face, a map $\varphi$ as in Proposition~\ref{easy:prop} gives an
$\{\alpha,\beta,\gamma,\delta\}$-coloring of the faces of $R$ with
the usual condition that adjacent faces should have different
colors. On the other hand, Lemma~\ref{Rfund:lem} implies the
converse, namely that any such coloring gives a map $\varphi$ as in
Proposition~\ref{easy:prop}, and hence a closed orientable
hyperbolic manifold. As an example, the classical L\"obell
manifold~\cite{Lo} is obtained from the polytope $R(6)$ described
above and shown in Fig.~\ref{Fig1}-left using the coloring shown on
the right in the same figure.

\begin{figure}[ht]
\begin{center}
\setlength{\unitlength}{0.4mm}
\begin{picture}(240,100)(0,10)
\special{em:linewidth 1.6pt} \put(-30,0){\begin{picture}(120,120)
\emline{40}{50}{1}{40}{70}{2} \emline{40}{70}{1}{60}{80}{2}
\emline{60}{80}{1}{80}{70}{2} \emline{80}{70}{1}{80}{50}{2}
\emline{80}{50}{1}{60}{40}{2} \emline{60}{40}{1}{40}{50}{2}
\emline{40}{50}{1}{30}{40}{2} \emline{40}{70}{1}{30}{80}{2}
\emline{60}{80}{1}{60}{90}{2} \emline{80}{70}{1}{90}{80}{2}
\emline{80}{50}{1}{90}{40}{2} \emline{60}{40}{1}{60}{30}{2}
\emline{40}{25}{1}{30}{40}{2} \emline{40}{25}{1}{60}{30}{2}
\emline{20}{60}{1}{30}{80}{2} \emline{20}{60}{1}{30}{40}{2}
\emline{40}{95}{1}{60}{90}{2} \emline{40}{95}{1}{30}{80}{2}
\emline{80}{95}{1}{60}{90}{2} \emline{80}{95}{1}{90}{80}{2}
\emline{100}{60}{1}{90}{80}{2} \emline{100}{60}{1}{90}{40}{2}
\emline{80}{25}{1}{60}{30}{2} \emline{80}{25}{1}{90}{40}{2}
\emline{30}{10}{1}{0}{60}{2} \emline{0}{60}{1}{30}{110}{2}
\emline{30}{110}{1}{90}{110}{2} \emline{90}{110}{1}{120}{60}{2}
\emline{120}{60}{1}{90}{10}{2} \emline{90}{10}{1}{30}{10}{2}
\emline{30}{10}{1}{40}{25}{2} \emline{0}{60}{1}{20}{60}{2}
\emline{30}{110}{1}{40}{95}{2} \emline{90}{110}{1}{80}{95}{2}
\emline{120}{60}{1}{100}{60}{2} \emline{90}{10}{1}{80}{25}{2}
\put(75,35){\makebox(0,0)[cc]{$1$}}
\put(90,60){\makebox(0,0)[cc]{$2$}}
\put(75,85){\makebox(0,0)[cc]{$3$}}
\put(45,85){\makebox(0,0)[cc]{$4$}}
\put(30,60){\makebox(0,0)[cc]{$5$}}
\put(45,35){\makebox(0,0)[cc]{$6$}}
\put(60,20){\makebox(0,0)[cc]{$7$}}
\put(100,40){\makebox(0,0)[cc]{$8$}}
\put(100,80){\makebox(0,0)[cc]{$9$}}
\put(60,100){\makebox(0,0)[cc]{$10$}}
\put(20,80){\makebox(0,0)[cc]{$11$}}
\put(20,40){\makebox(0,0)[cc]{$12$}}
\put(60,60){\makebox(0,0)[cc]{$13$}}
\put(110,100){\makebox(0,0)[cc]{$14$}}
\end{picture}}
\put(140,0){\begin{picture}(120,120) \emline{40}{50}{1}{40}{70}{2}
\emline{40}{70}{1}{60}{80}{2} \emline{60}{80}{1}{80}{70}{2}
\emline{80}{70}{1}{80}{50}{2} \emline{80}{50}{1}{60}{40}{2}
\emline{60}{40}{1}{40}{50}{2} \emline{40}{50}{1}{30}{40}{2}
\emline{40}{70}{1}{30}{80}{2} \emline{60}{80}{1}{60}{90}{2}
\emline{80}{70}{1}{90}{80}{2} \emline{80}{50}{1}{90}{40}{2}
\emline{60}{40}{1}{60}{30}{2} \emline{40}{25}{1}{30}{40}{2}
\emline{40}{25}{1}{60}{30}{2} \emline{20}{60}{1}{30}{80}{2}
\emline{20}{60}{1}{30}{40}{2} \emline{40}{95}{1}{60}{90}{2}
\emline{40}{95}{1}{30}{80}{2} \emline{80}{95}{1}{60}{90}{2}
\emline{80}{95}{1}{90}{80}{2} \emline{100}{60}{1}{90}{80}{2}
\emline{100}{60}{1}{90}{40}{2} \emline{80}{25}{1}{60}{30}{2}
\emline{80}{25}{1}{90}{40}{2} \emline{30}{10}{1}{0}{60}{2}
\emline{0}{60}{1}{30}{110}{2} \emline{30}{110}{1}{90}{110}{2}
\emline{90}{110}{1}{120}{60}{2} \emline{120}{60}{1}{90}{10}{2}
\emline{90}{10}{1}{30}{10}{2} \emline{30}{10}{1}{40}{25}{2}
\emline{0}{60}{1}{20}{60}{2} \emline{30}{110}{1}{40}{95}{2}
\emline{90}{110}{1}{80}{95}{2} \emline{120}{60}{1}{100}{60}{2}
\emline{90}{10}{1}{80}{25}{2}
\put(75,35){\makebox(0,0)[cc]{$\beta$}}
\put(90,60){\makebox(0,0)[cc]{$\gamma$}}
\put(75,85){\makebox(0,0)[cc]{$\delta$}}
\put(45,85){\makebox(0,0)[cc]{$\beta$}}
\put(30,60){\makebox(0,0)[cc]{$\gamma$}}
\put(45,35){\makebox(0,0)[cc]{$\delta$}}
\put(60,20){\makebox(0,0)[cc]{$\gamma$}}
\put(100,40){\makebox(0,0)[cc]{$\delta$}}
\put(100,80){\makebox(0,0)[cc]{$\beta$}}
\put(60,100){\makebox(0,0)[cc]{$\gamma$}}
\put(20,80){\makebox(0,0)[cc]{$\delta$}}
\put(20,40){\makebox(0,0)[cc]{$\beta$}}
\put(60,60){\makebox(0,0)[cc]{$\alpha$}}
\put(110,100){\makebox(0,0)[cc]{$\alpha$}}
\end{picture} }
\end{picture}
\end{center}
\caption{The polytope $R(6)$ and a coloring of its faces.}
\label{Fig1}
\end{figure}
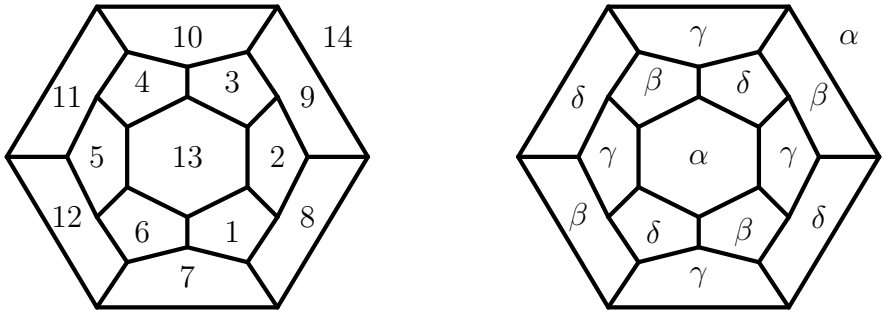

We generalize this example by considering for each $n\geqslant 5$
the right-angled polytope $R(n)$ in $\Hy^3$ with $(2n+2)$ faces, two
of which (viewed as the upper and  lower bases) are regular
$n$-gons, while the lateral surface is given by 2n pentagons,
arranged as one easily imagines. Note that $R(5)$ is the
right-angled dodecahedron. As in Fig.~\ref{Fig1} for $R(6)$, we
number the faces of $R(n)$ so that:

\begin{itemize}
\item The upper and lower bases have numbers $2n+1$ and $2n$
respectively; \item The pentagons adjacent to the upper basis are
cyclically numbered $1,\ldots,n$; \item The pentagons adjacent to
the lower basis are cyclically numbered $n+1,\ldots,2n$ in the same
verse as the previous ones, with pentagon $n+1$ adjacent to
pentagons $1$ and $n$.
\end{itemize}

Now define $g_i\in{\rm Isom}(\Hy^3)$ as the reflection in the
plane containing the $i$-th face of $R(n)$, and $G(n)$ as the
group generated $\{g_i\}_{i=1}^{2n+2}$. According to
Lemma~\ref{Rfund:lem} a presentation of $G(n)$ is obtained by
adding the relations:
$$\begin{array}{ll}
g_i^2 &
i=1,\ldots,2n+2  \\
\left[g_{2n+1},g_i\right] &
i=1,\ldots,n \\
\left[g_{2n+2},g_{n+i}\right]\quad\  &
i=1,\ldots,n \\
\left[g_i,g_{i+1}\right]
& i=1,\ldots, 2n - 1 \\
\left[g_1,g_n\right] & \\
\left[g_i,g_{n+i}\right] &
i=1,\ldots,n \\
\left[g_{n+1},g_{2n}\right] & \\
\left[g_i,g_{n+1+i}\right] &
i=1,\dots,n-1.
\end{array}$$

We now define the class 
{\em L\"obell manifolds} of order $n$ as follows:
$$\begin{array}{rl}
\calL(n)=\Big\{\Hy^3/\Ker(\varphi): &
\varphi:G(n)\to (\Z/_{\!2})^3\ {\rm epimorphism},\\
& \Ker(\varphi)<{\rm Isom}^+(\Hy^3),\\
& \Ker(\varphi)\ \hbox{is\ torsion-free}\Big\}.
\end{array}$$
Each element of $\calL(n)$ is a closed orientable hyperbolic
3-manifold with volume equal to 8 times the volume of $R(n)$.
According to the above discussion (and the 4-color theorem!) the set
$\calL(n)$ is non-empty for all $n\geqslant 5$. The classical
L\"obell manifold constructed in~\cite{Lo} and described above
belongs to $\calL(6)$. We now provide an upper bound for the
complexity of the elements of $\calL(n)$.

\begin{lemma} For all $n\geqslant 5$ and $M\in\calL(n)$ we have
$c(M) \leqslant 32(2n -1)$.
\end{lemma}

\begin{proof}
By construction $M$ is built by gluing together in pairs the faces
of $8$ copies of $R(n)$. Each face is a $k$-gon with $k\in\{5,n\}$
and we can subdivide it into $k-2$ triangles by inserting $k-3$
diagonals, in such a way that the gluing respects the subdivision.
The boundary of each of the 8 copies of $R(n)$ is now subdivided
into $2(n-2)+2n(5-2)=4(2n-1)$ triangles. Taking the cone from the
center, we can subdivide the copy of $R(n)$ itself into the same
number of tetrahedra, which yields the desired upper bound
immediately.
\end{proof}

To give lower complexity estimates we will employ the hyperbolic
volume. Let us denote by $\ell_n$ the common value of ${\rm vol}(M)$
as $M$ varies in $\calL(n)$, and recall the definition of the {\em
Lobachevskii function}
$$\Lambda(x) = - \int\limits_0^x \log | 2 \sin (t) | \, {\rm d} t.$$
The following was established in~\cite{Ve2}:

\begin{thm}\label{theoremL3.3}
For all $n\geqslant 5$ we have
$$
\ell_n =  4 n \left( 2 \Lambda (\theta)  +  \Lambda \left( \theta +
\frac{\pi}{n} \right)  +  \Lambda \left( \theta - \frac{\pi}{n}
\right)  -  \Lambda \left( 2 \theta - \frac{\pi}{2} \right) \right)
,
$$
where $$\theta  =  \frac{\pi}{2}  -  \arccos  \left(\frac{1}{2 \cos
( \pi / n )} \right).$$
\end{thm}

This theorem implies that the volume of the classical L\"obell
manifold is equal to $48.184368\ldots $\ . In addition, it allows us
to determine the asymptotic behavior of $\ell_n$ as $n$ tends to
infinity. Indeed, in the limit the angle $\theta$ of the statement
tends to $\pi/6$, and using the fact that $\Lambda$ is continuous
and odd we have:

\begin{cor} \label{corollaryL4.1}
As $n$ tends to $\infty$ we have
$$\ell_n  \sim  10n\cdot  v_3$$
where $v_3=2 \Lambda (\pi / 6 ) = 1.014\dots$ is the volume of the
regular ideal tetrahedron in $\Hy^3$.
\end{cor}

To apply this result we establish the following general fact:

\begin{prop}\label{obvious:prop}
If $M$ is a closed orientable hyperbolic manifold then $${\rm
vol}(M)<c(M)\cdot v_3.$$
\end{prop}

\begin{proof}
Denote $c(M)$ by $k$. Since $M$ is irreducible and not one of the
exceptional manifolds $\mathbb S^3$, $\R\matP^3$, and $L(3,1)$,
there exists a realization of $M$ as a gluing of $k$ tetrahedra.
Denoting by $\Delta$ the abstract tetrahedron, this realization
induces continuous maps $\sigma_i:\Delta\to M$ for $i=1,\ldots,k$
describing how the tetrahedra appear in $M$ after the gluing. Note
that each $\sigma_i$ is injective on the interior of $\Delta$ but
perhaps not on the boundary. Since the gluings used to pair the
faces of the tetrahedra in the construction of $M$ are simplicial,
we see that $\sum_{i=1}^k \sigma_i$ is a singular $3$-cycle, which
of course represents the fundamental class $[M]\in H_3(M;\Z)$.

Let us now consider the universal covering $\Hy^3\to M$. Since
$\Delta$ is simply connected, we can lift $\sigma_i$ to a map
$\widetilde{\sigma}_i:\Delta\to\Hy^3$. We denote now by
$\widetilde{\tau}_i:\Delta\to\Hy^3$ the simplicial map which agrees
with $\widetilde{\sigma}_i$ on the vertices, where geodesic convex
combinations are used in $\Hy^3$ to define the notion of
`simplicial.' Let $\tau_i:\Delta\to M$ by the composition of
$\widetilde{\tau}_i$ with the projection $\Hy^3\to M$. One can
easily see that $\sum_{i=1}^k \tau_i$ is again a singular $3$-cycle
in $M$. Using this fact and taking convex combinations in $\Hy^3$
one can actually see that the cycles $\sum_{i=1}^k \sigma_i$ and
$\sum_{i=1}^k \tau_i$ are homotopic. In particular, since the first
cycle represents $[M]$, the latter also does, which implies that
$\bigcup_{i=1}^k\tau_i(\Delta)$ is equal to $M$, otherwise
$\sum_{i=1}^k \tau_i$ would be homotopic to a map with 2-dimensional
image.

Let us now note that $\widetilde{\tau}_i(\Delta)$ is a compact
geodesic tetrahedron in $\Hy^3$, therefore its volume is less than
$v_3$, see~\cite{lectures}. Moreover the volume of $\tau_i(\Delta)$
is at most equal to the volume of $\widetilde{\tau}_i(\Delta)$,
because the projection $\Hy^3\to M$ is a local isometry, and the
volume of $M$ is at most the sum of the volumes of the
$\tau_i(\Delta)$'s, because we have shown above that $M$ is covered
by the $\tau_i(\Delta)$'s. This concludes the proof.
\end{proof}

\begin{cor}
For sufficiently large $n$ and $M\in\calL(n)$ we have
 $c(M)>10 n$.
\end{cor}

The following theorem summarizes the complexity bounds we have
obtained in the present section:

\begin{thm}
For sufficiently large $n$ and $M\in\calL(n)$ the complexity of $M$
satisfies the following inequalities:
$$
10 n  < c(M)  \leqslant  32(2n -1) .
$$
\end{thm}

\section{Fibonacci manifolds}

In this section we consider the compact orientable hyperbolic
3-manifolds whose fundamental groups are the Fibonacci groups,
introduced by J.~Conway in~\cite{Con1}. There is one such group
$F(2,n)$ for each $n\geqslant 3$, and a presentation of it is given
by
$$
F(2,n)= \langle x_1,x_2,\dots ,x_n :\ x_i x_{i+1}x_{i+2}^{-1},\
i=1,\ldots,n \rangle
$$
where indices are understood modulo $n$. It was shown in~\cite{HKM1}
that for each $n\geqslant 4$ 
the group $F(2,2n)$ is isomorphic to a discrete cocompact
subgroup of ${\rm Isom}^+(\Hy^3)$, the group of
orientation-preserving isometries of hyperbolic 3-space. We will
need below to refer to several details of the construction
of~\cite{HKM1}, so we recall it here.

\begin{figure}[ht]
\begin{center}
\unitlength=0.5mm
\begin{picture}(200,160)(0,30)
\thicklines \put(10,80){\line(1,0){160}}
\put(30,120){\line(1,0){160}} \put(10,80){\line(1,2){20}}
\put(50,80){\line(1,2){20}} \put(90,80){\line(1,2){20}}
\put(130,80){\line(1,2){20}} \put(170,80){\line(1,2){20}}
\put(50,80){\line(-1,2){20}} \put(90,80){\line(-1,2){20}}
\put(130,80){\line(-1,2){20}} \put(170,80){\line(-1,2){20}}
\put(110,20){\line(-5,3){100}} \put(110,20){\line(-1,1){60}}
\put(110,20){\line(-1,3){20}} \put(110,20){\line(1,3){20}}
\put(110,20){\line(1,1){60}} \put(110,180){\line(-4,-3){80}}
\put(110,180){\line(-2,-3){40}} \put(110,180){\line(0,-3){60}}
\put(110,180){\line(2,-3){40}} \put(110,180){\line(4,-3){80}}
\put(125,20){\makebox(0,0)[cc]{$R$}}
\put(125,185){\makebox(0,0)[cc]{$Q$}}
\put(30,95){\makebox(0,0)[cc]{$F_5^{*}$}}
\put(70,95){\makebox(0,0)[cc]{$F_7^{*}$}}
\put(110,95){\makebox(0,0)[cc]{$F_1^{*}$}}
\put(150,95){\makebox(0,0)[cc]{$F_3^{*}$}}
\put(50,105){\makebox(0,0)[cc]{$F_6^{*}$}}
\put(90,105){\makebox(0,0)[cc]{$F_8^{*}$}}
\put(130,105){\makebox(0,0)[cc]{$F_2^{*}$}}
\put(170,105){\makebox(0,0)[cc]{$F_4^{*}$}}
\put(60,55){\makebox(0,0)[cc]{$F_6$}}
\put(85,55){\makebox(0,0)[cc]{$F_8$}}
\put(110,55){\makebox(0,0)[cc]{$F_2$}}
\put(135,55){\makebox(0,0)[cc]{$F_4$}}
\put(75,145){\makebox(0,0)[cc]{$F_7$}}
\put(100,145){\makebox(0,0)[cc]{$F_1$}}
\put(120,145){\makebox(0,0)[cc]{$F_3$}}
\put(145,145){\makebox(0,0)[cc]{$F_5$}}
\put(25,125){\makebox(0,0)[cc]{$P_8$}}
\put(65,125){\makebox(0,0)[cc]{$P_2$}}
\put(100,125){\makebox(0,0)[cc]{$P_4$}}
\put(155,125){\makebox(0,0)[cc]{$P_6$}}
\put(195,125){\makebox(0,0)[cc]{$P_8$}}
\put(5,75){\makebox(0,0)[cc]{$P_7$}}
\put(45,75){\makebox(0,0)[cc]{$P_1$}}
\put(85,75){\makebox(0,0)[cc]{$P_3$}}
\put(135,75){\makebox(0,0)[cc]{$P_5$}}
\put(175,75){\makebox(0,0)[cc]{$P_7$}}
\end{picture}
\end{center}
\caption{The polytope $Y(4)$.} \label{Fig2}
\end{figure}
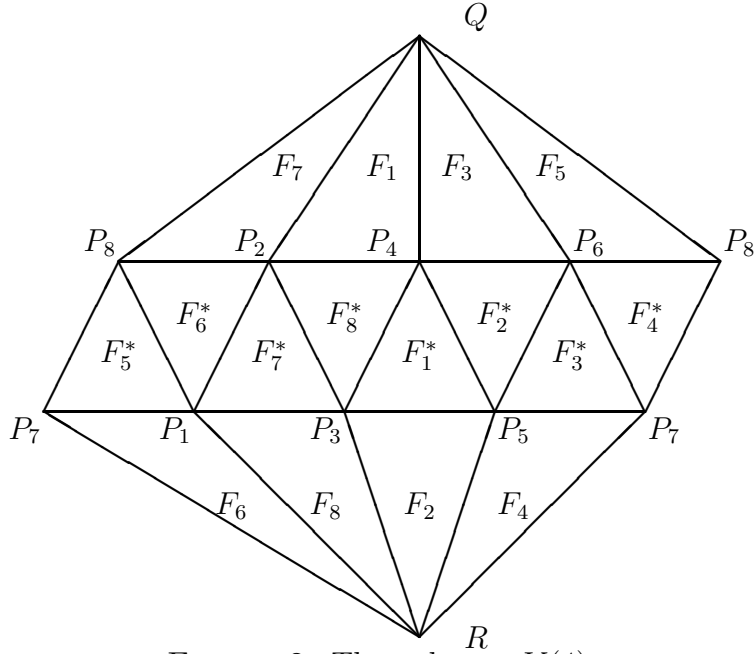

We fix $n \geqslant 4$ and we first define the order-$n$ antiprism
$A(n)$ as the polytope whose boundary is constructed as follows:
\begin{itemize}
\item Take $2n$ triangles and two polygons with $n$ faces;
\item Attach a different triangle to each edge of each of the two $n$-polygons;
\item Glue together the objects thus obtained by matching the free edges of the triangles
(there are two circles consisting of $n$ edges to glue together, so
there is essentially only one way to do so).

\end{itemize}
Now we define $Y(n)$ as the polytope obtained from $A(n)$ by
attaching an $n$-pyramid to each of the bases. In particular, $Y(5)$
is the icosahedron. We remark that in general $Y(n)$ has $2n+2$
vertices, $6n$ edges, and $4n$ triangular faces, and we denote the
vertices by $Q$, $R$, $P_1, \dots , P_{2n}$ and the faces by $F_1,
\dots , F_n$, and $F_1^{*}, \dots , F_n^{*}$, as shown in
Fig.~\ref{Fig2} for $n=4$.

We define now a face-pairing on $Y(n)$ under which each face $F_i$
is glued to the face $F_i^*$ via a simplicial homeomorphism $s_i :
F_i
 \longrightarrow  F_i^{*}$. We specify the action of $s_i$ by
 describing its action on the vertices. Namely, for odd $i$ we choose $s_i$
 so that
$$s_i  :  Q P_{i+1} P_{i+3}
 \longrightarrow P_{i+2} P_{i+3} P_{i+4}$$
whereas for even $i$ we choose it so that
$$s_i  :  R P_{i+1} P_{i+3}  \longrightarrow P_{i+2}
P_{i+3} P_{i+4}.$$ Note that if we choose an orientation of $Y(n)$
and orient the $F_i$ and $F_i^*$ accordingly, all the $s_i$'s are
orientation-reversing homeomorphisms. This
implies that the quotient of $Y(n)$ under the face-pairing is a
manifold except perhaps at the vertices, and that the projection
restricted to each open edge of $Y(n)$ is injective. In particular,
we can describe how the various edges of $Y(n)$ are cyclically
arranged around an edge in the quotient. These cycles are actually
easy to describe: for odd $i$ we have
\begin{equation}
Q P_{i+1}  \stackrel{s_i}{\longrightarrow}  P_{i+2} P_{i+3}
\stackrel{ s_{i-1}^{-1}}{\longrightarrow}  P_{i} P_{i+2}
\stackrel{s_{i-2}^{-1}}{ \longrightarrow}  Q P_{i+1} , \label{fm3}
\end{equation}
and for even $i$ we have 
\begin{equation}
R P_{i+1}  \stackrel{s_i}{\longrightarrow}  P_{i+2} P_{i+3}
\stackrel{ s_{i-1}^{-1}}{\longrightarrow}  P_{i} P_{i+2}
\stackrel{s_{i-2}^{-1}}{ \longrightarrow}  R P_{i+1} . \label{fm4}
\end{equation}

It was shown in~\cite{HKM1} that $Y(n)$ can be realized in a unique
way (up to isometry) as a compact polytope in hyperbolic space
$\Hy^3$ in such a way that:
\begin{itemize}
\item Each of the faces of $Y(n)$ is an equilateral triangle;
\item The sums of the dihedral angles
corresponding to the cycles (\ref{fm3}) and (\ref{fm4}) is
 $2 \pi$;
\item $Y(n)$ has a cyclic symmetry of order $n$ with axis $QR$ and
an orientation-reversing symmetry which permutes $Q$ and  $R$.
\end{itemize}
We will henceforth identify $Y(n)$ with such a realization in
$\Hy^3$. Since all the faces of $Y(n)$ are congruent, each
face-pairing $s_i$ can be realized in a unique fashion as an
orientation-preserving isometry of $\Hy^3$, and we will denote this
isometry also by $s_i$. The condition that the total dihedral angle
around the edge-cycles (\ref{fm3}) and (\ref{fm4}) is $2\pi$ easily
implies that
\begin{equation} s_i  s_{i+1}  =  s_{i+2},\qquad i = 1, \dots , 2n
\label{fm5}
\end{equation}
where indices are understood modulo $2n$. More precisely,
Poincar\'e's polyhedron theorem~\cite{EpPet,Vin} implies that:
\begin{itemize}
\item The subgroup of
${\rm Isom}^+(\Hy^3)$ generated by the $s_i$'s is isomorphic to
$F(2,2n)$;
\item This group is discrete and torsion-free, and $Y(n)$ is a fundamental domain for
its action on $\Hy^3$;
\item The
quotient of $\Hy^3$ under this action is a $3$-manifold.
\end{itemize}
We will denote from now on by $M(n)$ the closed hyperbolic
3-manifold thus obtained, and call it the $n$-th {\em Fibonacci
manifold}. It was remarked in~\cite{HLM1} that $M(n)$ is the
$n$-fold cyclic covering of $\mathbb S^3$ branched over the
figure-eight knot $4_1$.

\begin{lemma} {\em For $n \geqslant 4$ we have
$c(M(n)) \leqslant 3n$.}
\end{lemma}

\begin{proof}
For each triangular face of $Y(n)$ not containing $Q$ we can
construct the tetrahedron with vertex at $Q$ and basis at that face.
This gives a decomposition of $Y(n)$ into $3n$ tetrahedra, whence a
singular triangulation of $M(n)$ with $3n$ tetrahedra, whence the
conclusion at once.\end{proof}

To estimate the complexity of $M(n)$ from below we use the following
formula for its volume established in~\cite{VeMe1}:

\begin{thm}\label{theoremF1}
For $n \geqslant 4$ we have
$$
{\rm vol}(M(n)) = 2 n ( \Lambda(a_n + b_n) + \Lambda(a_n - b_n) ),
$$
where $b_n = {\pi}/{n}$ and $a_n = (1/2)\cdot \arccos(\cos(2a_n) -
{1}/{2})$.
\end{thm}

This result allows us to determine the asymptotic behavior of the
volume of the Fibonacci manifold $M(n)$. Indeed, as $n$ tends to
infinity, we see that $b_n$ converges to $0$ and $a_n$ converges to
$\pi / 6$. Using the fact that $\Lambda$ is continuous we deduce the
following:

\begin{cor}
As $n$ tends to $\infty$ we have
$${\rm vol}(M(n))  \sim  2n\cdot  v_3$$
where $v_3=2 \Lambda (\pi / 6 ) = 1.014\dots$ is the volume of the
regular ideal tetrahedron in $\Hy^3$.
\end{cor}

This result together with Proposition~\ref{obvious:prop} implies:

\begin{cor}
For sufficiently large $n$ we have
 $c(M(n))>2n$.
\end{cor}

The following theorem summarizes the complexity bounds we have
obtained in the present section:

\begin{thm}
For sufficiently large $n$ the complexity of the Fibonacci manifold
$M(n)$ satisfies the following inequalities:
$$
2 n  < c(M(n))  \leqslant  3n .
$$
\end{thm}

\end{document}